\documentclass[11pt]{amsart}
\usepackage{amsthm,amsmath,amssymb}
\usepackage{mathtools}
\usepackage{graphicx}
\usepackage{caption}
\graphicspath{{Figure/}}
\numberwithin{equation}{section}
\numberwithin{figure}{section}
\theoremstyle{plain}

\newtheorem{thm}{Theorem}[section]
\newtheorem{lem}[thm]{Lemma}

\theoremstyle{definition}
\theoremstyle{remark}

\begin{document}
\setlength{\abovedisplayskip}{10pt}
\setlength{\belowdisplayskip}{10pt}

\title{Lozenge tilings of hexagons with holes on three crossing lines}

\author{Seok Hyun Byun}
\address{Department of Mathematics, Indiana University, Bloomington}
\email{byunse@indiana.edu}

\maketitle
\begin{abstract}
The enumeration of lozenge tilings of hexagons with holes has received much attention during the last three decades. One notable feature is that a lot of the recent development involved Kuo's graphical condensation. Motivated by Ciucu, Lai and Rohatgi's work on tilings of hexagons with a removed triad of bowties, in this paper, we show that the ratio of numbers of lozenge tilings of two more general regions is expressed as a simple product formula. Our proof does not involve the graphical condensation method. The proof is short and direct. We also provide a corresponding formula for cyclically symmetric lozenge tilings. Several previous results can be easily deduced from our generalization.
\end{abstract}

\section{Introduction}
David and Tomei's bijection [15] between plane partitions fitting inside a box and lozenge tilings of a corresponding hexagon on the triangular lattice allows one to interpret MacMahon's classical theorem [30] as follows: The number of lozenge tilings of a hexagon with side lengths \(a\), \(b\), \(c\), \(a\), \(b\), \(c\) (clockwise from top) is given by the following beautiful product formula:

\begin{equation}
    \prod_{i=1}^{\ a}\prod_{j=1}^{\ b}\prod_{k=1}^{\ c} \frac{i+j+k-1}{i+j+k-2}.
\end{equation}

In the 1990s, several people generalized this result by enumerating the number of lozenge tilings of a hexagon with a triangular hole at the center (see Ciucu [2], Ciucu, Eisenk\"{o}lbl, Krattenthaler and Zare [7], Gessel and Helfgott [17], Okada and Krattenthaler [31]).

Later, further generalizations have been discovered. Some of them generalized it by increasing the number of triangular holes (see Ciucu [3]) or changing the shape of the hole at the center (see Ciucu [4], Ciucu and Krattenthaler [10], Lai and Rohatgi [27]). Others generalized it by putting some holes along the boundary of the hexagon (see Ciucu and Lai [11], Lai [21, 22]) or holes at both center and boundary (see Lai [23]). Figure 1.1 shows just a handful from the dozens of regions on the triangular lattice that can be found in the literature whose number of lozenge tilings is given by a simple product formula. The first picture on the top corresponds to formula (1.1). The second picture on the top is an example of \textit{cored hexagons}, regions introduced by Ciucu, Eisenk\"{o}lbl, Krattenthaler and Zare that form the subject of [7]. The third picture on the top is a special example\footnote{This picture is a special case because, in the original context [3], Ciucu allowed hexagons to have multiple triangles removed.} of \textit{symmetric hexagons with certain triangles along its symmetry axis removed}, regions studied by Ciucu in [3]. Cored hexagons mentioned earlier were then extended in [10] by Ciucu and Krattenthaler to hexagons with a so-called \textit{shamrock} removed from its center (see the first picture in the middle). The second picture in the middle is an instance of a family of regions, \textit{a hexagon with three dents}, introduced by Lai in [21]. The last picture in the middle is an example of \textit{axial shamrocks} that were studied independently by Ciucu [6] and Lai and Rohatgi in [27]. Cored hexagons were once more generalized by Ciucu to hexagons with a so-called \textit{fern} removed from the center (see [4] and the first picture on the bottom). The second picture on the bottom represents an example of \textit{doubly intruded hexagons}, regions studied by Ciucu and Lai in [11]. Later, Lai generalized two previous results and investigated \textit{hexagons with three ferns removed} (see [23] and the last picture on the bottom).

\begin{figure}
    \centering
    \includegraphics[width=11.5cm]{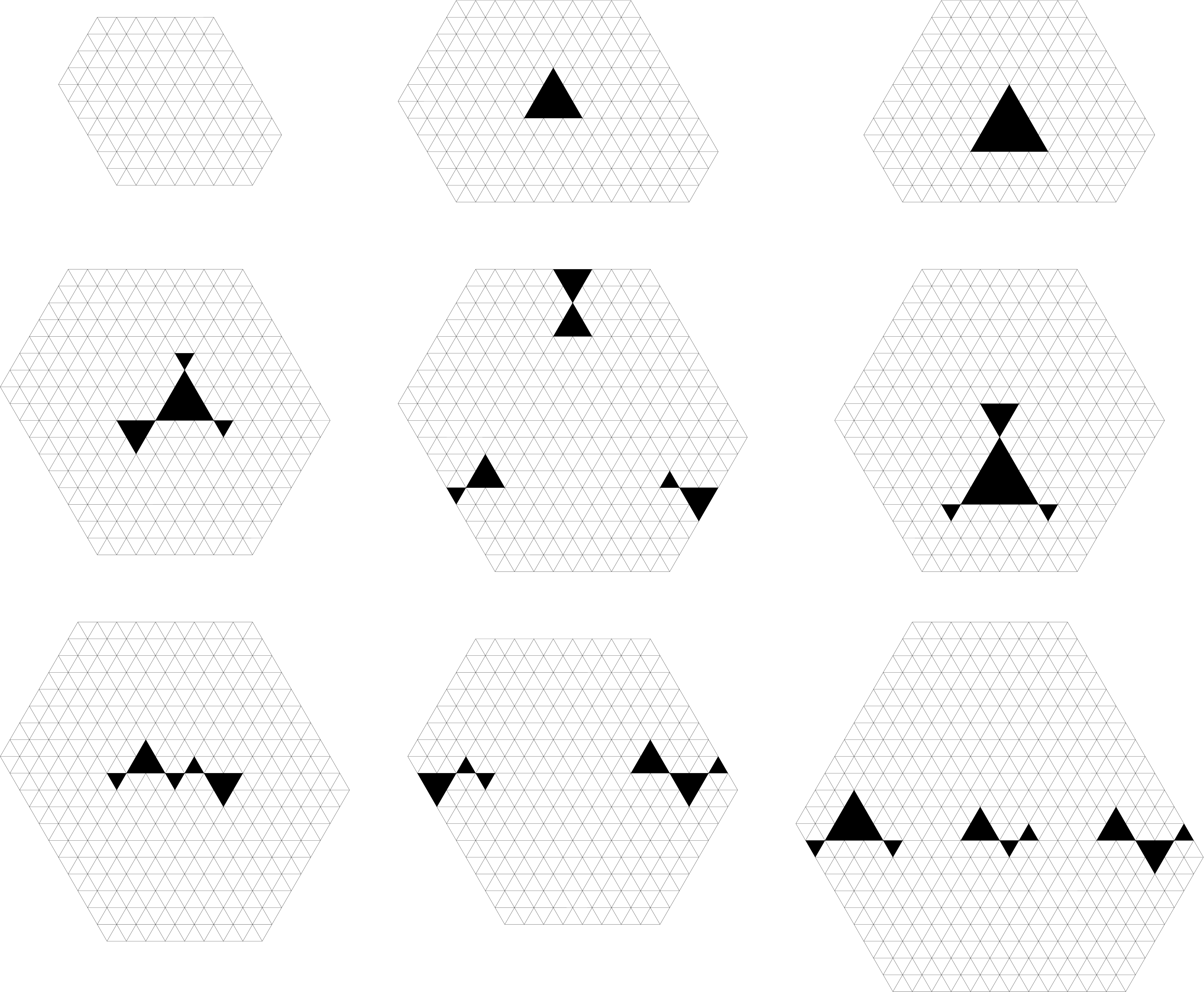}
    \caption{Various regions whose number of lozenge tilings are given by simple product formulas.}
\end{figure}

More recently, simple product formulas for the ratio of the number of lozenge tilings of two related regions were also found (see Ciucu, Lai and Rohatgi [12], Condon [13, 14], Lai [24, 25, 26], and Lai and Rohatgi~[28]). In particular, Lai and Rohatgi's \textit{shuffling theorem} [28] and Ciucu, Lai and Rohatgi's \textit{bowtie squeezing theorem} [12] provided a simple proof of many lozenge tilings enumeration problems mentioned in the previous paragraph. For example, in Figure 1.1, lozenge tilings enumeration of three figures in the middle could be easily obtained by bowtie squeezing theorem, and shuffling theorem provided a simple proof of that of three figures on the bottom.

One notable feature is that a lot of the recent development involves Kuo's graphical condensation method [19, 20]. It is a powerful tool in the field of enumeration of tilings. If one can guess the formula for the number of tilings of a family of regions and if this family of regions is general enough, this method could allow proving the formula by induction. However, simple product formulas from recent results call for more direct and straightforward proofs of them. 

Motivated by Ciucu, Lai and Rohatgi's work [12] on tilings of hexagons with a removed triad of bowties (in which proofs are based on Kuo condensation), in this paper, we show that the ratio of the numbers of lozenge tilings of two more general regions is expressed as a simple product formula. We will also see how this identity can give unified proofs (or explanations) of some recent results proven by various arguments involving Kuo's graphical condensation method. Our simple argument enables us to relate this result to the enumeration of cyclically symmetric lozenge tilings (i.e., tilings invariant under rotation by 120\(^{\circ}\)). We will in fact provide a corresponding formula for cyclically symmetric lozenge tilings for the same kind of regions. It provides in particular a simple proof for the enumeration of cyclically symmetric lozenge tilings of a hexagon with a \textit{shamrock} (a certain 4-lobed structure) removed from the center. This was first proved by Ciucu in [5]; that proof uses Ciucu and Fischer's work [8], which involves the graphical condensation method. By contrast, the current paper does not use graphical condensation in its proofs. Our arguments and Kuo condensation complement each other, as explained at the end of Section 2.

This paper is organized as follows. In Section 2, we state the main theorem, give a geometric interpretation and see how it provides unified proofs for some previous results from the literature that were originally proved in different ways. In Section 3, we state a key lemma, we prove it and then show how to use it to prove the main theorem. We also provide formulas for two more symmetry classes at the end of Section 3. 

\section{Statement of Main Results}

In this paper, we consider bounded regions on a triangular lattice. Without loss of generality, we draw the lattice so that one family of the lines is horizontal. A \textit{lozenge} is the union of two adjacent unit triangles on the lattice. Given a region, a \textit{lozenge tiling} of the region is a collection of lozenges that covers it without gaps or overlaps. There are three types of lozenges that one can consider: Left-, vertical- and right-lozenges (see Figure 2.1). We now describe the family of regions that we will deal with in this paper.

For non-negative integers \(n\) and \(x\), we consider the hexagon whose side lengths are \(n\), \(n+x\), \(n\), \(n+x\), \(n\), \(n+x\) (clockwise from top). The three long diagonals determine an up-pointing triangle of side length \(x\) at the center of the hexagon. We remove the triangle of size \(x\) from the hexagon and denote the remaining region by \(H_{n,x}\). Let \(U\), \(L\) and \(R\) be the upper, bottom left and bottom-right vertex of this triangular hole, respectively. Label the unit segments on the line from \(U\) to the top-right vertex of the hexagon by~\(1,2,\dotsc,n\) from inside out. We use the analogous labeling on the five lines connecting \(R\) and the right vertex of the hexagon, \(R\) and bottom-right vertex of the hexagon, \(L\) and bottom-left vertex of the hexagon, \(L\) and left vertex of the hexagon and \(U\) and top-left vertex of the hexagon (see Figure~2.2).

\begin{figure}
    \centering
    \includegraphics[width=7cm]{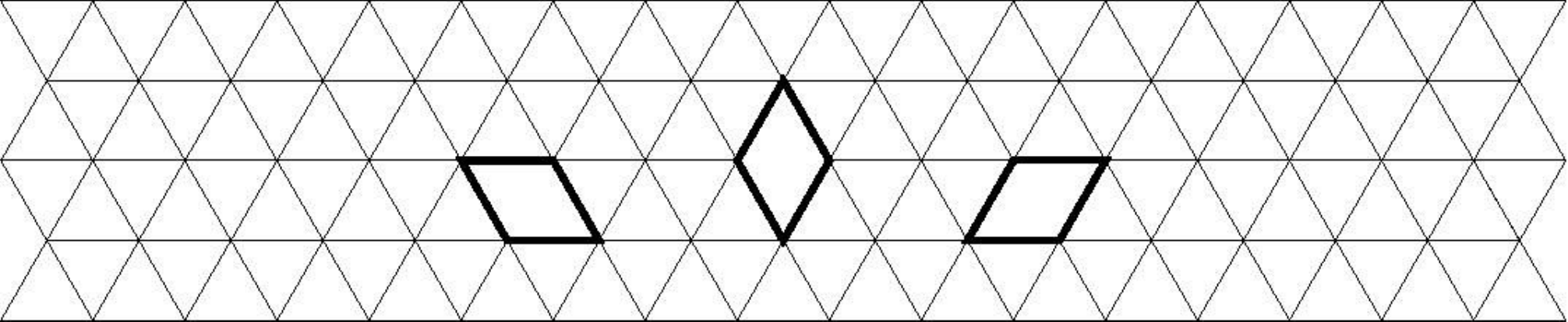}
    \caption{Left-, vertical- and right-lozenge on triangular lattice}
\end{figure}

\begin{figure}
    \centering
    \includegraphics[width=8.8cm]{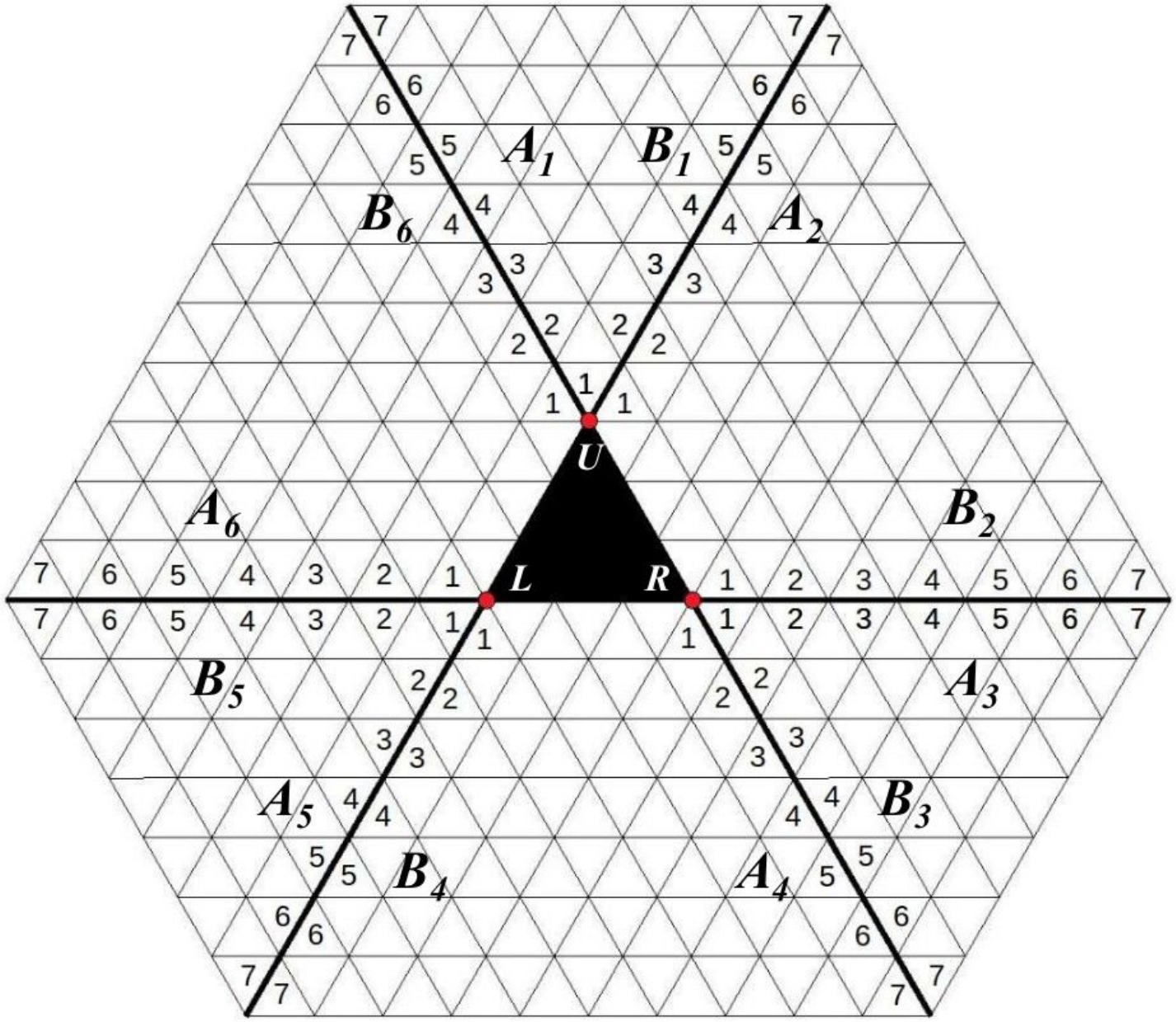}
    \caption{\(H_{7,3}\); \(U\), \(L\), \(R\), labels of unit segments and portion each set \(A_1,\dotsc, B_6\) occupies are represented.}
\end{figure}

Note that these 6 lines decompose \(H_{n,x}\) into 6 subregions: 3 trapezoids and 3 triangles. Let \(\textbf{A}\coloneqq (A_1, A_2,\dotsc, A_6)\) and \(\textbf{B}\coloneqq (B_1, B_2,\dotsc, B_6)\), where \(A_1\), \(A_2,\dotsc\), \(A_6\), \(B_1\), \(B_2,\dotsc\), \(B_6\) are subsets of \([n]\coloneqq\{1,2,\dotsc,n\}\). The region \(H_{n,x}(\textbf{A},\textbf{B})\) in the theorem below is obtained from \(H_{n,x}\) by removing a collection of unit triangles from along the sides of these subregions as follows. Consider \(A_1\) and \(B_1\) as sets of labels on the left and right sides of the triangle on top (see Figure 2.2). Similarly, consider the remaining \(A_i\)'s and \(B_i\)'s as sets of labels as indicated in Figure 2.2. 

\begin{figure}
    \centering
    \includegraphics[width=10.5cm]{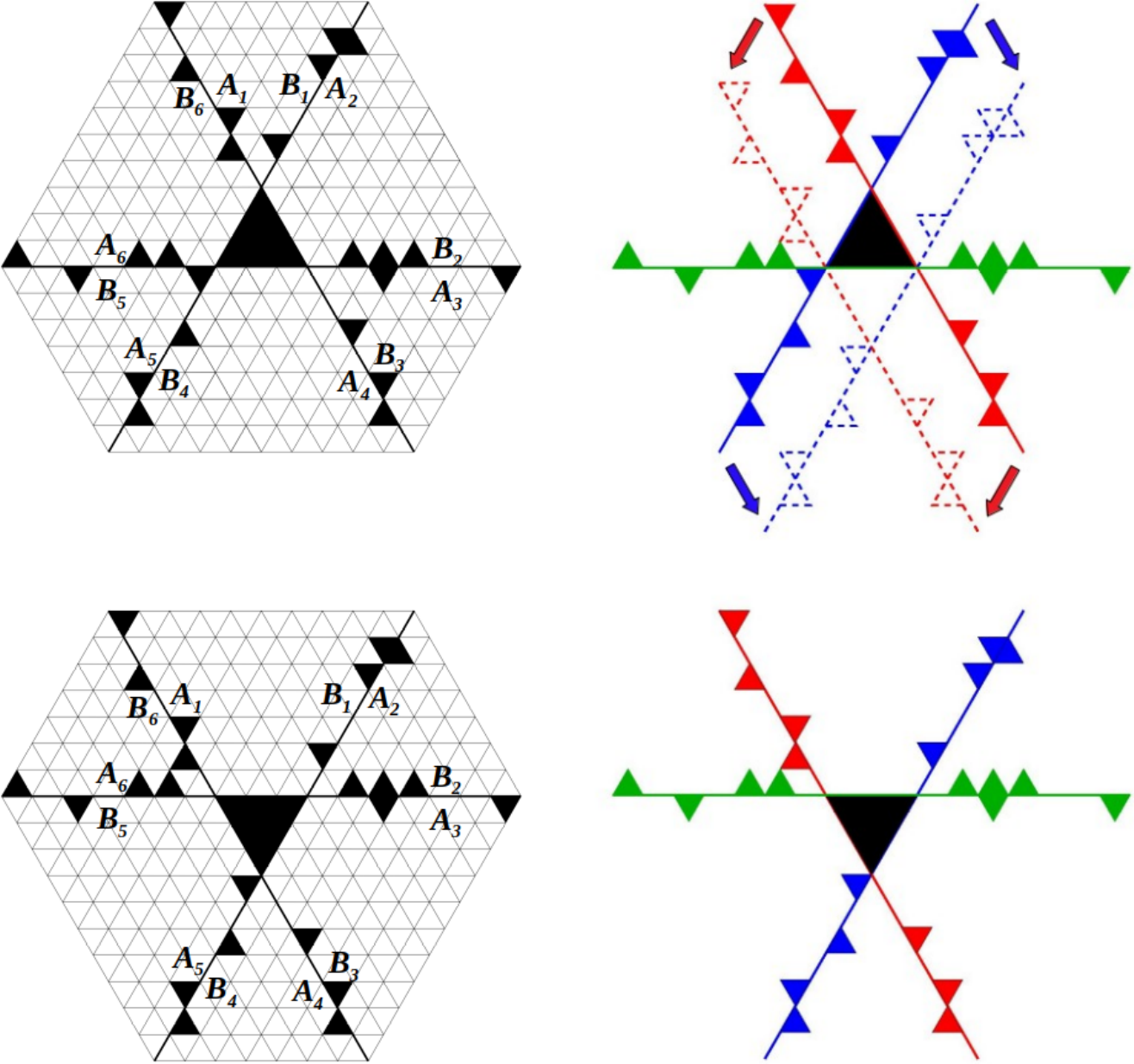}
    \caption{The figure \(H_{7,3}(\textbf{A},\textbf{B})\) with \(A_1=\{3,7\}\), \(B_1=\{2,5,6\}\), \(A_2=\{6\}\), \(B_2=\{2,3,4\}\), \(A_3=\{3,7\}\), \(B_3=\{3,5\}\), \(A_4=\{6\}\), \(B_4=\{3,6\}\), \(A_5=\{1,5\}\), \(B_5=\{5\}\), \(A_6=\{2,3,7\}\) and \(B_6=\{2,5\}\) is on the top left. Follow the figures clockwise from top left to see the snowflake flipping process. The figure on the bottom left is \(\overline{H}_{7,3}(\textbf{A},\textbf{B})\)}
\end{figure}

We define \(H_{n,x}(\textbf{A},\textbf{B})\) to be the region obtained from \(H_{n,x}\) by removing the unit triangles touching the labeled sides at positions specified by \(A_1, A_2,\dotsc, A_6, B_1, B_2,\dotsc, B_6\). Given the resemblance to a snowflake, we group the unit holes specified by the \(A_i\)'s and \(B_i\)'s into three \textit{dendrites}: The \textit{horizontal dendrite} consists of the horizontal long diagonal and the unit holes specified by the sets \(A_3\), \(A_6\), \(B_2\) and \(B_5\); the \textit{positive dendrite} consists of the positive slope long diagonal and the unit holes specified by the sets  \(A_2\), \(A_5\), \(B_1\) and \(B_4\); and the \textit{negative dendrite} consists of the negative slope long diagonal and the unit holes specified by the sets \(A_1\), \(A_4\), \(B_3\) and \(B_6\). We call the region \(H_{n,x}(\textbf{A},\textbf{B})\) a \textit{snowflake region}. Through out this paper, we will assume\footnote{If this assumption is missing, then there is a possibility that a single unit triangular hole located at the intersection of two dendrites is indexed by two sets \(A_{i}\) and \(B_{i}\) for some \(i\) either on \(H_{n,x}(\textbf{A},\textbf{B})\) or on \(\overline{H}_{n,x}(\textbf{A},\textbf{B})\) (The region \(\overline{H}_{n,x}(\textbf{A},\textbf{B})\) is defined shortly when we introduce snowflake flipping operation). Besides, one can easily check that either \(H_{n,x}(\textbf{A},\textbf{B})\) or \(\overline{H}_{n,x}(\textbf{A},\textbf{B})\) does not have any lozenge tilings if we miss this assumption.} that \(1\notin A_{i}\cap B_{i}\) for \(i=1,2,\dotsc,6\).

We now define an operation called \textit{snowflake flipping}, which transforms a given snowflake region into a new region (which becomes a new snowflake region when rotated by 180\(^{\circ}\)).

Notice that the region \(H_{n,x}(\textbf{A},\textbf{B})\) is determined by its dendrites (because the boundary hexagon is the convex hull of the union of the three dendrites). In order to visualize our definition of snowflake flipping, it will help to consider the unit triangular holes to be fastened to the axis of their corresponding dendrite. 

Our snowflake flipping operation is this: While keeping the horizontal dendrite fixed, translate the positive dendrite \(x\) units to the southeast along the lattice and the negative dendrite \(x\) units to the southwest along the lattice. Then, instead of the up-pointing triangle of side length \(x\), the newly placed dendrites determine a \textit{down}-pointing triangle, of the same side length \(x\). This triangle will be a hole in our new region. We enclose the dendrites in the convex hull of their union, which is a hexagon of side lengths \(n+x\), \(n\), \(n+x\), \(n\), \(n+x\), \(n\) (clockwise from top). We denote this new region by \(\overline{H}_{n,x}(\textbf{A},\textbf{B})\) (see Figure 2.3). We say \(\overline{H}_{n,x}(\textbf{A},\textbf{B})\) is obtained from \(H_{n,x}(\textbf{A},\textbf{B})\) by snowflake flipping.

Let \(A_{o}\) and \(A_{e}\) be multisets obtained by listing all elements of the three sets \(A_1\), \(A_3\), \(A_5\) and \(A_2\), \(A_4\), \(A_6\), respectively. Similarly, let \(B_{o}\) (resp., \(B_{e}\)) be multisets obtained by listing all elements of the three sets \(B_1\), \(B_3\), \(B_5\) (resp., \(B_2\), \(B_4\), \(B_6\)). Recall that the Pochhammer symbol \((a)_k\) is defined by \((a)_0\coloneqq 1\) and \((a)_k\coloneqq \displaystyle \prod_{i=0}^{k-1}(\alpha+i)\) for a positive integer \(k\).

For a region \(R\) on the triangular lattice, let \(M(R)\) be the number of its lozenge tilings. We say that the region \(R\) is \textit{cyclically symmetric} if the region is invariant under rotation by 120\(^{\circ}\). Our regions \(H_{n,x}(\textbf{A},\textbf{B})\) and \(\overline{H}_{n,x}(\textbf{A},\textbf{B})\) are cyclically symmetric if and only if \(A_1=A_3=A_5\), \(A_2=A_4=A_6\), \(B_1=B_3=B_5\) and \(B_2=B_4=B_6\) hold. A lozenge tiling of a cyclically symmetric region is called \textit{cyclically symmetric} if the tiling is invariant under rotation by 120\(^{\circ}\). For such a region \(R\), let \(M_r(R)\) be the number of cyclically symmetric lozenge tilings of it. The main result of this paper is the following.

\begin{thm}
Let \(n\) and \(x\) be non-negative integers. Suppose \(A_1,\) \(\dotsc,\) \(A_6,\) \(B_1,\) \(\dotsc,\) \(B_6\subseteq[n]\) satisfy \(1\notin A_{i}\cap B_{i}\), for \(i=1,2,\dotsc,6\).

\emph{(a)}. If \(H_{n,x}(\textbf{A},\textbf{B})\) has a lozenge tiling, then
\begin{equation}
    \frac{M(\overline{H}_{n,x}(\textbf{A},\textbf{B}))}{M(H_{n,x}(\textbf{A},\textbf{B}))}=\frac{\displaystyle \prod_{a\in A_{o}}(a)_x \ \displaystyle \prod_{b\in B_{o}}(b)_x}{\displaystyle \prod_{a\in A_{e}}(a)_x \ \displaystyle \prod_{b\in B_{e}}(b)_x}.
\end{equation}

\emph{(b)}. If \(H_{n,x}(\textbf{A},\textbf{B})\) is cyclically symmetric and has a cyclically symmetric lozenge tiling, then
\begin{equation}
    \frac{M_{r}(\overline{H}_{n,x}(\textbf{A},\textbf{B}))}{M_{r}(H_{n,x}(\textbf{A},\textbf{B}))}=\frac{\displaystyle \prod_{a\in A_{1}}(a)_x \ \displaystyle \prod_{b\in B_{1}}(b)_x}{\displaystyle \prod_{a\in A_{2}}(a)_x \ \displaystyle \prod_{b\in B_{2}}(b)_x}=\sqrt[3]{\frac{M(\overline{H}_{n,x}(\textbf{A},\textbf{B}))}{M(H_{n,x}(\textbf{A},\textbf{B}))}}.
\end{equation}
\end{thm}

\textbf{Geometric interpretation.} Note that the unit holes corresponding to the sets \(A_{1}\), \(A_{3}\), \(A_{5}\), \(B_{1}\), \(B_{3}\) and \(B_{5}\) are contained in the three triangular subregions of \(H_{n,x}(\textbf{A},\textbf{B})\) (see Figure 2.2), and they are all down-pointing. On the other hand, the unit holes corresponding to the remaining six sets are contained in the three trapezoidal subregions of \(H_{n,x}(\textbf{A},\textbf{B})\), and they are all up-pointing.

Based on this, we can give the following geometric interpretation to the right-hand side of equation (2.1). 
Define the distance  \(d(A,B))\) between two unit triangles \(A\) and \(B\) supported on a common lattice line \(\ell\) to be the Euclidean distance (with unit being the side length of a unit triangle) between the centers of their sides that are along \(\ell\). For example, in Figure~2.4, \(d(A,B))\)=\(5\) and \(d(B,C))\)=\(3\).

Let \(C\) be a triangle of side length \(k\) supported\footnote{We say that a lattice triangle \(T\) is supported on the lattice line \(\ell\) if \(\ell\) contains one side of \(T\).} on the lattice line \(l\). Define the \textit{projection of \(C\) onto \(l\)} to be the set of \(k\) unit triangles inside \(C\) that have a side along  \(l\) (see Figure 2.5).

In our theorem, for a unit hole labeled by \(a\), the corresponding factor \((a)_{x}=a(a+1)\dotsc(a+x-1)\) in (2.1) represents the product of distances between that unit hole and the unit triangles in the projection of the central hole of size \(x\) onto the lattice line supporting these two holes. 

When a unit triangle \(A\) and a triangle of any size \(C\) are supported on the same lattice line, define \(C_A\) to be the projection of \(C\) onto that line. Then, equation (2.1) can be stated as

\begin{figure}
    \centering
    \includegraphics[width=11cm]{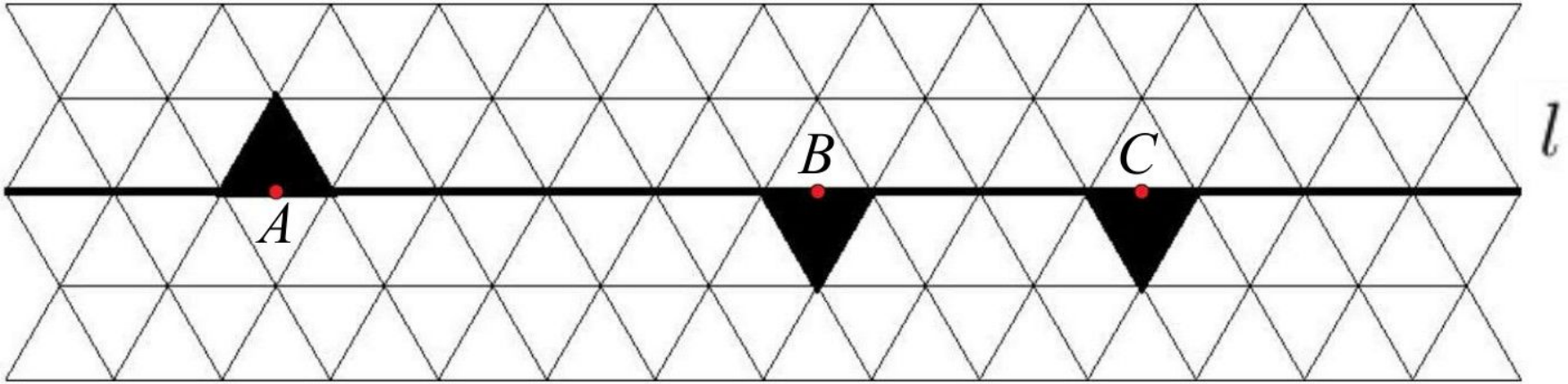}
    \caption{Unit triangles hanging on the same line \(l\). Red dots represent midpoints of unit segments that each triangle shares with the line.}
\end{figure}
\begin{figure}
    \centering
    \includegraphics[width=11cm]{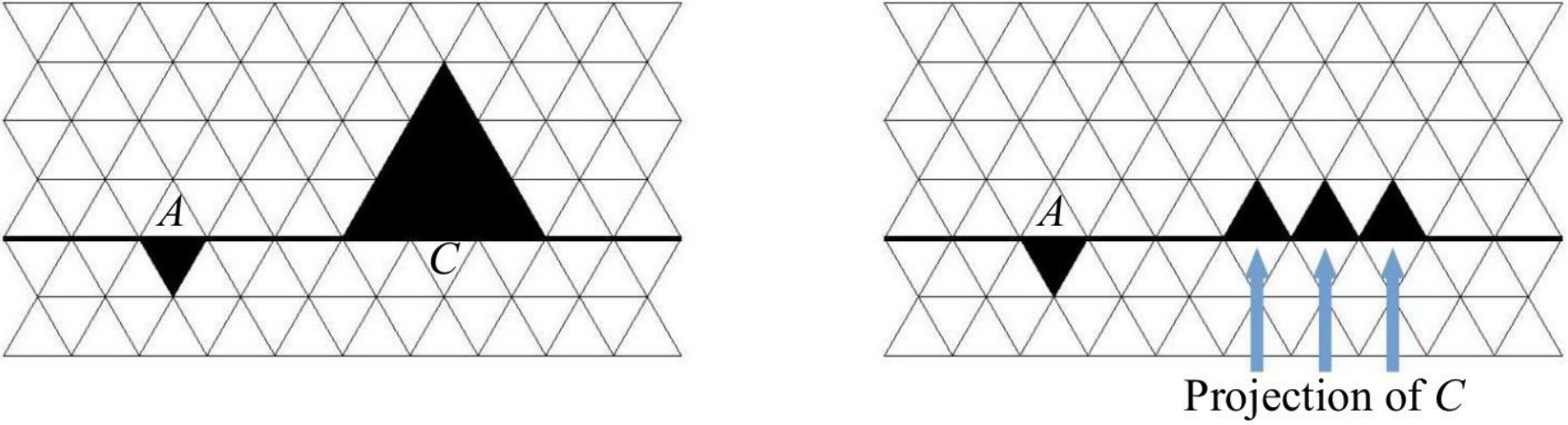}
    \caption{Two triangles \(A\) and \(C\) on the same line (left) and figure after projecting \(C\) onto the line (right).}
\end{figure}

\begin{equation}
    \frac{M(\overline{H}_{n,x}(\textbf{A},\textbf{B}))}{M(H_{n,x}(\textbf{A},\textbf{B}))}=\frac{\displaystyle \prod_{\blacktriangledown}\prod_{\triangle\in X_{\blacktriangledown}} d(\blacktriangledown, \triangle)}{\displaystyle \prod_{\blacktriangle}\prod_{\triangle\in X_{\blacktriangle}} d(\blacktriangle, \triangle)}
\end{equation}
where \(X\) represents the triangular hole of size \(x\) at the center in the theorem, the outer product at numerator runs over all down-pointing unit holes hanging on the dendrites, and similarly for outer product on the denominator.

This form brings out the geometric meaning of the formula on the right hand side of (2.1).

\smallskip
For a region \(R\) on the triangular lattice, a \textit{forced lozenge} is a lozenge-shaped subregion of \(R\) that is always covered by a single lozenge in all tilings of the region. Thus, if we make a lozenge-shaped hole on such a spot, the number of lozenge tilings of the new region is equal to that of the original region. In the figures of this paper, some of the force lozenges are indicated by a shading.

\begin{figure}
    \centering
    \includegraphics[width=11cm]{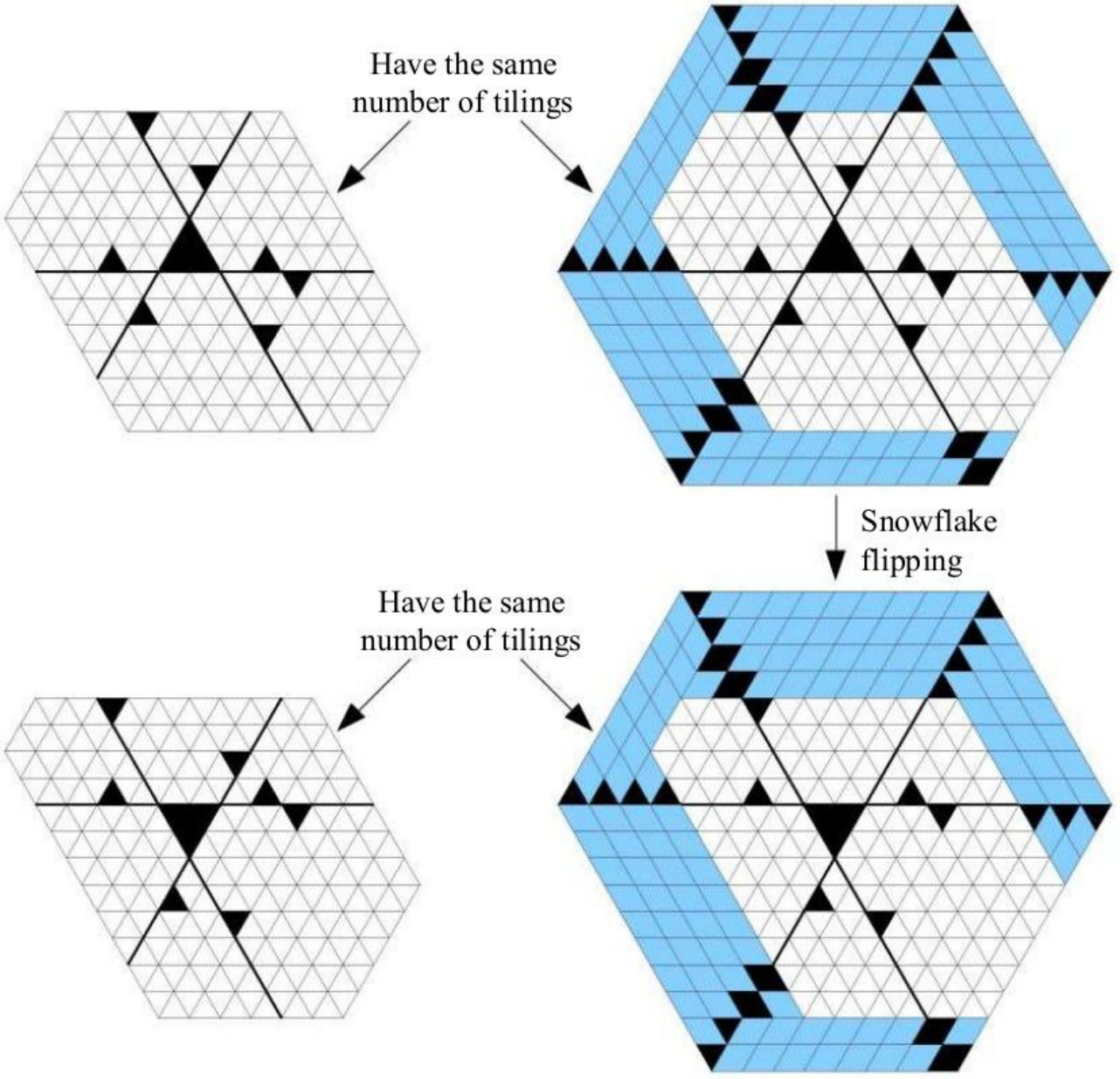}
    \caption{Snowflake flipping operation for arbitrary hexagon with holes on three crossing lines in general positions. On the two right figures, forced lozenges generated by newly added holes are represented by shading; read the figures clockwise from top left.}
\end{figure}

One may think that the shape of the hexagon and the positions of the three dendrites (namely a cyclically symmetric hexagon and its three long diagonals) in Theorem 2.1 are special. However, by an idea used previously by Ciucu and Lai [11], it turns out that our theorem provides in fact a product formula for any hexagon and any choice of three crossing lines. Indeed, we can see from Figure 2.6, we can always find a cyclically symmetric hexagon with holes on the three long diagonals that has the same number of lozenge tilings as the original region (compare the two figures on top). This can be done by extending the three dendrite axes and placing additional holes on their extensions as indicated in the figure. If we apply the Theorem 2.1 to the top right region, we have a simple formula for the ratio between the number of lozenge tilings of the two regions on the right. However, if we discard some forced lozenges from the bottom right region, we get a region whose central triangle is flipped compared to that of the original region (compare the two figures on the left). Since the ratio of the number of lozenge tilings of the two figures on the left is the same as that of the two figures on the right, the theorem provides a formula for the ratio of the number of tilings of the two regions on the left. This argument also gives a natural explanation for the result stated as Corollary 2.2 in Rosengren's paper [32]: The corollary was the particular case when there is only one hole in the hexagon (and no other unit holes along the dendrites).

We now indicate how one can use the main theorem to deduce three results which were proved originally by separate arguments. The first two of them were originally proved by using Kuo's graphical condensation method.

\begin{figure}
    \centering
    \includegraphics[width=11cm]{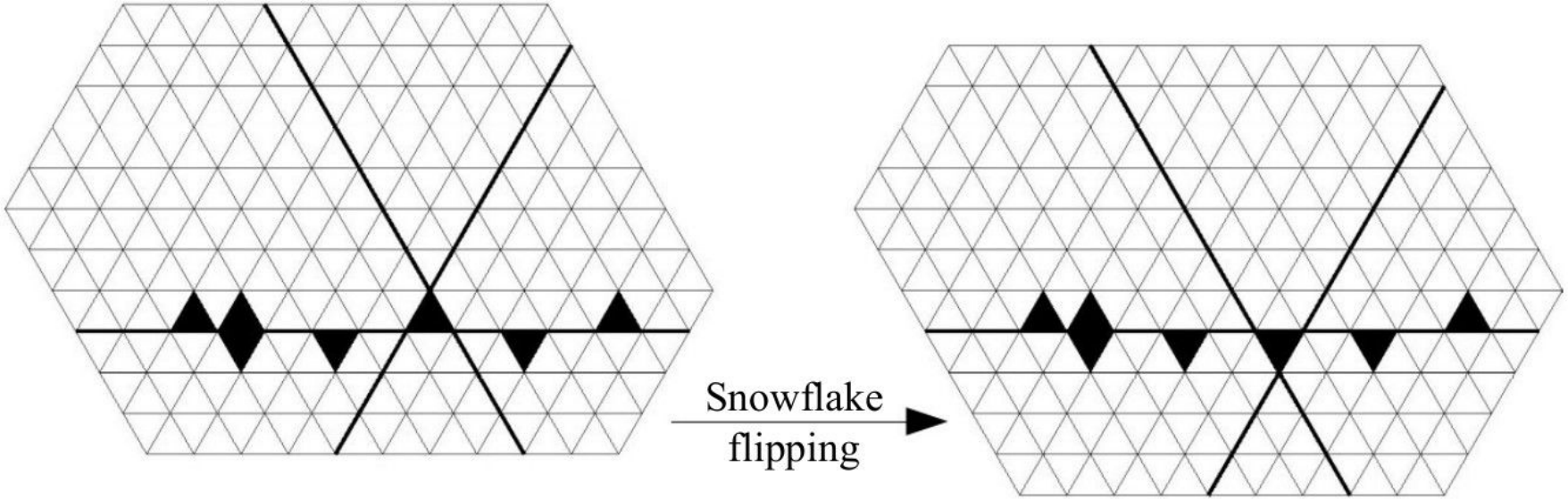}
    \caption{Flipping a unit triangle}
\end{figure}

\begin{figure}
    \centering
    \includegraphics[width=11cm]{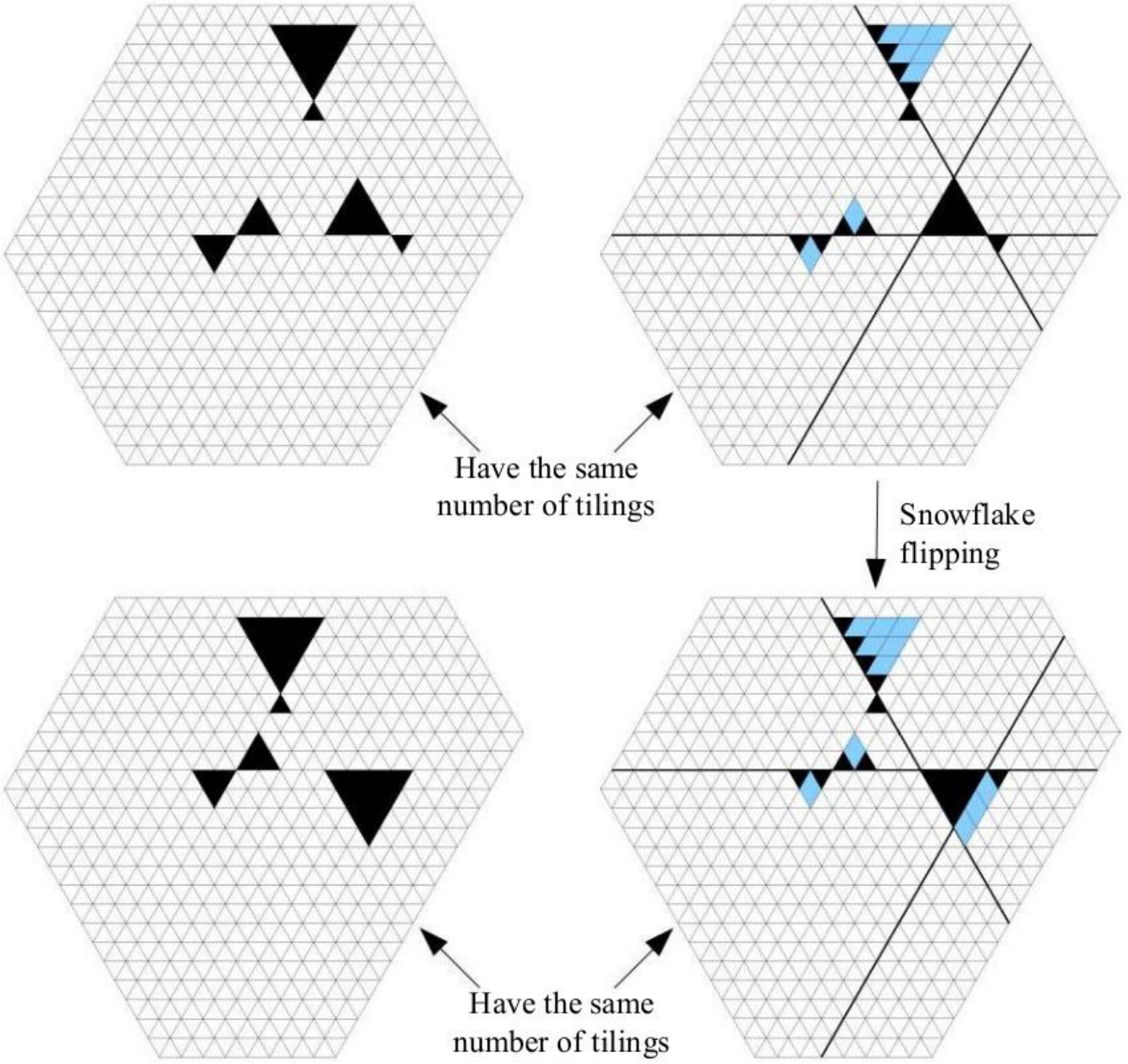}
    \caption{Sqeezing a bowtie}
\end{figure}

\begin{figure}
    \centering
    \includegraphics[width=11cm]{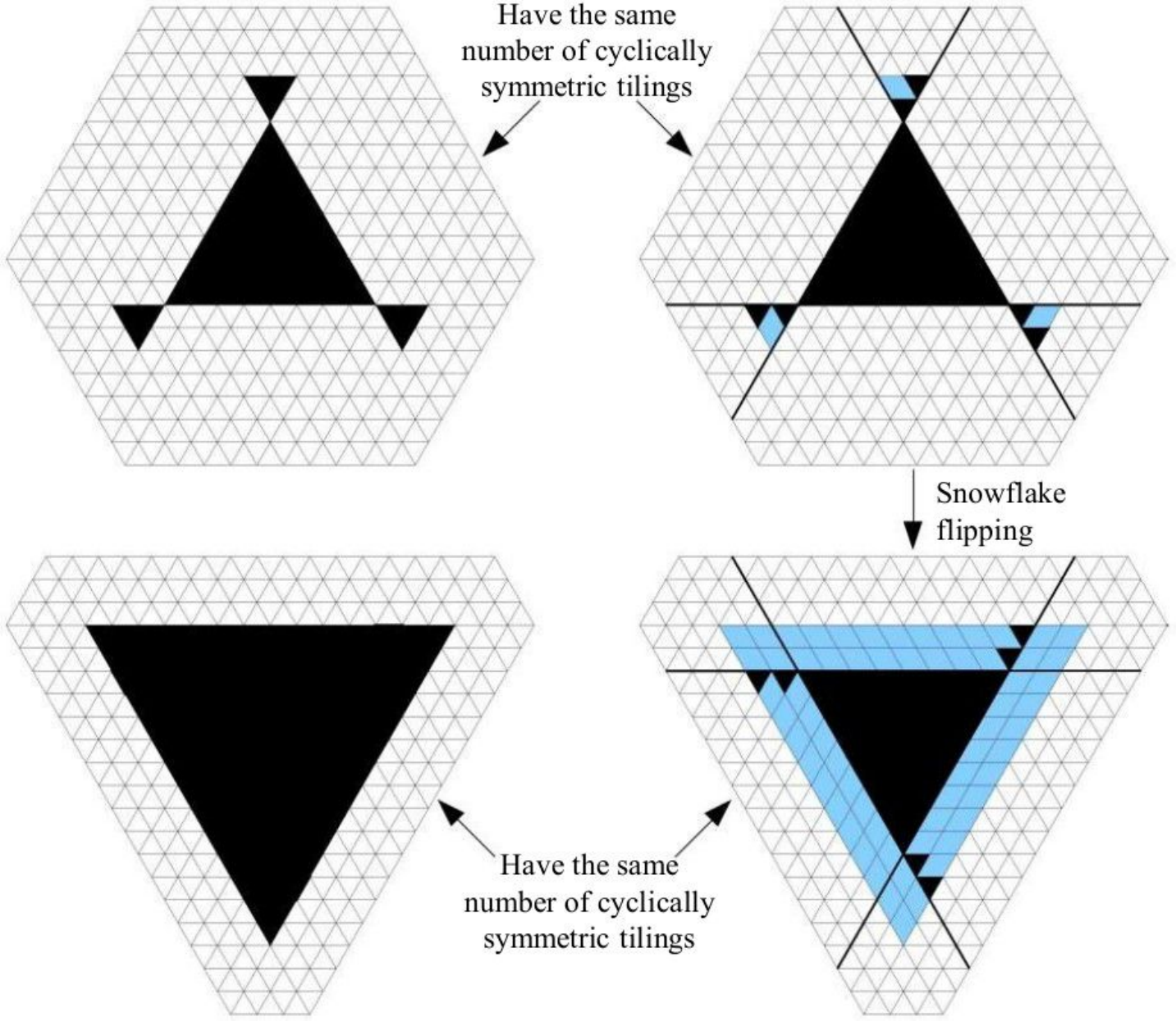}
    \caption{Application of snowflake flipping to a cyclically symmetric hexagon with shamrock removed at the center}
\end{figure}

The first is Lai and Rohatgi's shuffling theorem of lozenge tilings (see~[28]). This was already proved in a different way, independently, by the author [1] and Fulmek [16]. It also follows from a particular case of our main theorem, when the size of the central hole is \(1\), and there are no other holes on the positive and negative dendrites. Indeed, as we can see\footnote{We need to see both Figure 2.7 and Figure 2.6 to understand this argument. In fact, the left picture in Figure 2.7 is not a snowflake region. However, the argument that we used in Figure 2.6 enables us to apply the main theorem in this case. Later, we also view Figure 2.8 and Figure 2.9 together with Figure 2.6.} from Figures~2.7 and 2.6, Theorem 2.1 gives the ratio between the number of tilings of a shuffled region with a single unit triangle flipped and that of the original region. The case of an arbitrary shuffled region is then obtained by applying such a unit triangle flipping repeatedly.

The second is Ciucu, Lai and Rohatgi's bowtie squeezing theorem (see~[12]), which was in fact the original motivation for the current paper. Any hexagon with a removed triad of bowties (see the top left picture in Figure 2.8 for an example) can be viewed as being obtained from a convenient snowflake region after removing some forced lozenges (see Figures 2.8 and 2.6). By applying Theorem 2.1 and putting lozenge-shaped holes on some forced lozenges (indicated by a shading in the bottom right picture in Figure 2.8), we obtain an expression for the ratio between the number of tilings of the region with one bowtie completely squeezed out (bottom left picture in Figure 2.8) and the original region. The general case follows by applying such single bowtie squeezing two more times.

The third corollary of Theorem 2.1 concerns the enumeration of cyclically symmetric lozenge tilings of a hexagon with a shamrock (see Figure 2.9 and [10]) removed from the center, which was worked out by Ciucu in [5]. This can be deduced 
from previous work of Ciucu and Krattenthaler~[9] and Theorem 2.1(b). In [9], the authors enumerated cyclically symmetric lozenge tilings of a hexagon with a single triangular hole at the center (their proof is based on the non-intersecting lattice path approach and does not involve Kuo condensation). After removing forced lozenges, a hexagon with a shamrock removed from the center can be viewed as a snowflake region (see Figures 2.9 and 2.6). Thus, part (b) of Theorem 2.1 provides a formula for the ratio between the number of cyclically symmetric lozenge tilings of a hexagon with a single triangular hole and that of a hexagon with a shamrock hole (see Figure 2.9). Ciucu’s result follows then by combining this with [9].

Our theorem and Kuo's graphical condensation method complement each other in the following sense. While our theorem provides a simple product formula for the ratio between numbers of lozenge tilings of two regions, it does not seem to be able to prove a formula for a single region. On the other hand, the graphical condensation method can often be used to prove such a formula. For instance, in [10], Ciucu and Krattenthaler provided a simple proof for the number of lozenge tilings of a hexagon with a shamrock removed from the center by using the graphical condensation method. (A special case of their result is when a triangle is removed from the center; see~[7]).

\section{Proof of the main theorem}

Let \(L_{n,x}\) be the boundary parallelogram in the left picture in Figure 3.1, where the horizontal sides have length \(n\) and the oblique sides have length \(n+x\). Split \(L_{n,x}\) into a triangle and a trapezoid, as shown in the figure. The region \(L_{n,x}(P,Q,R,S)\) in the lemma below is obtained from \(L_{n,x}\) by removing the collection of unit triangles from along the sides of this triangle and trapezoid as follows. Label the unit segments along the left and right sides of the triangle from bottom to top by \(1,2,\dotsc,n\); let \(P, Q \subseteq [n]\), and consider them as sets of labels on these two sides, respectively. Similarly, label the unit segments along the top left and bottom sides of the trapezoid from left to right by \(1,2,\dotsc,n\); let \(R, S \subseteq [n]\) denote sets of labels along these sides, respectively. We define \(L_{n,x}(P,Q,R,S)\) to be the region obtained from \(L_{n,x}\) by removing the unit triangles touching the labeled sides at positions specified by \(P,Q,R,S\). The region \(\overline{L}_{n,x}(P,Q,R,S)\) is defined analogously using the picture on the right in Figure 3.1. 

Our proof of Theorem 2.1 is based on the following result.

\begin{figure}
    \centering
    \includegraphics[width=11cm]{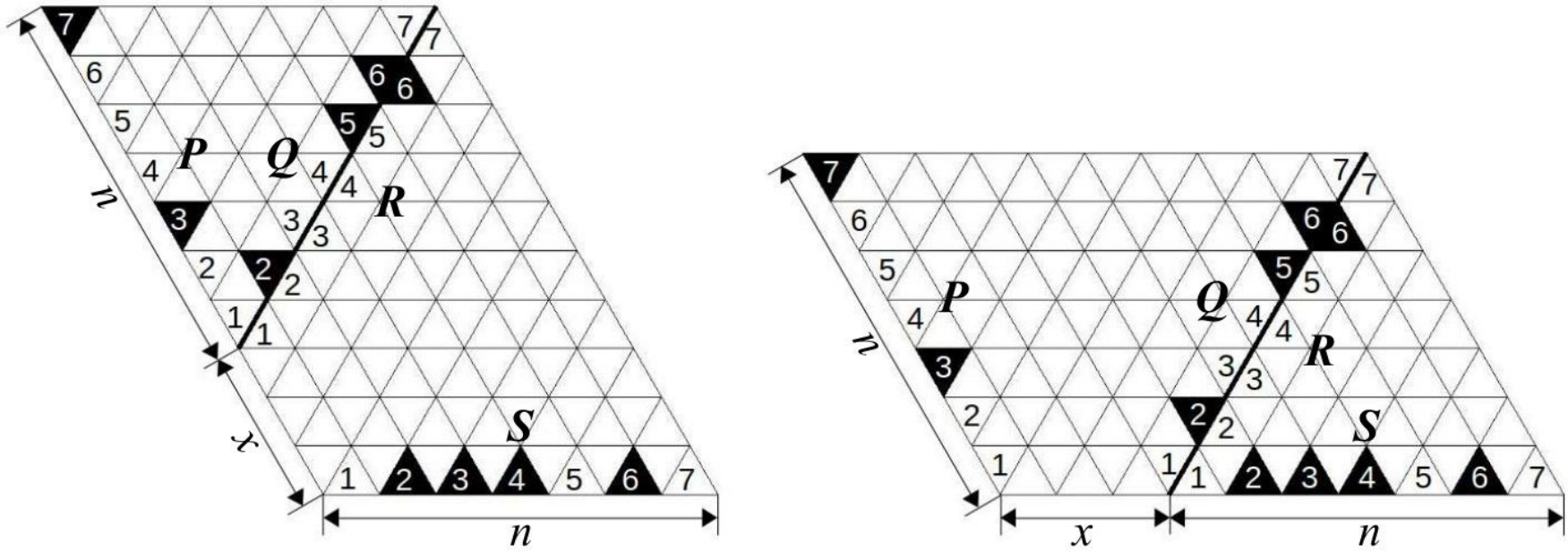}
    \caption{\(L_{7,3}(P,Q,R,S)\) (left) and \(\overline{L}_{7,3}(P,Q,R,S)\) (right) with \(P=\{3,7\}\), \(Q=\{2,5,6\}\), \(R=\{6\}\) and \(S=\{2,3,4,6\}\)}
\end{figure}

\begin{lem} 
Let \(n\) and \(x\) be non-negative integers. Suppose \(P\), \(Q\), \(R\) and \(S \subseteq [n]\) satisfy \(1 \notin P\cap Q\) and \(1 \notin R\cap S\). Then if the region \(L_{n,x}(P,Q,R,S)\) can be tiled by lozenges, we have

\begin{equation}
    \frac{M(\overline{L}_{n,x}(P,Q,R,S))}{M(L_{n,x}(P,Q,R,S))}=\frac{\displaystyle \prod_{p\in P}(p)_x \ \prod_{q\in Q}(q)_x}{\displaystyle \prod_{r\in R}(r)_x \ \prod_{s\in S}(s)_x}.
\end{equation}
\end{lem}

Note that the unit holes labeled by \(P\) and \(Q\) point down, while those labeled by \(R\) and \(S\) point up.

The idea of the proof is the following. We use the well known interpretation of lozenge tilings as families of non-intersecting lattice paths on \(\mathbb{Z}^2\) (see Figure 3.2). For both our \(L\)- and \(\overline{L}\)-regions, partition the collection of such families of non-intersecting lattice paths according to the points on the line \(y=x\) that they pass through. Since the segments supporting the holes indexed by \(Q\) and \(R\) in \(L_{n,x}(P,Q,R,S)\) and \(\overline{L}_{n,x}(P,Q,R,S)\) look exactly the same, the classes of these two partitions are naturally paired up. The key element is that the ratio between the cardinalities of paired partition classes turns out (using the Lindstr\"{o}m-Gessel-Viennot theorem [18, 29]) to be equal to the right-hand side of (3.1), and does not depend on the choice of partition classes.

\medskip
\textit{Proof of Lemma 3.1.} We use the well-known interpretation of lozenge tiling as families of paths of lozenges, which in turn can be identified with families of non-intersecting lattice paths on \(\mathbb{Z}^2\) allowed to take steps in two cardinal directions (for us these will be south and east; see Figure 3.2 for an illustration. For now, the double points on the first bisector should be ignored). 

\begin{figure}
    \centering
    \includegraphics[width=11cm]{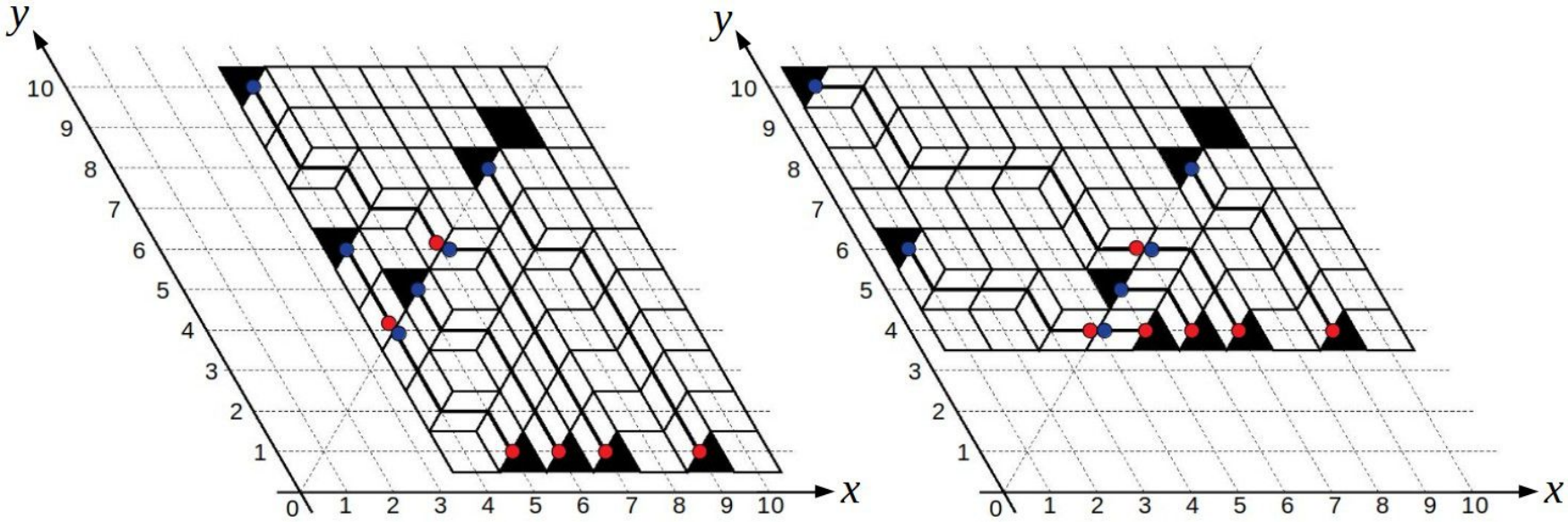}
    \caption{Non-intersecting lattice path interpretation of lozenge tilings of \(L_{n,x}(P,Q,R,S)\) (left) and \(\overline{L}_{n,x}(P,Q,R,S)\) (right), where parameters are same as in Figure 3.1 and \(U=\{1, 3\}\) (see the proof for meaning of the set \(U\)). Blue points represent starting points of lattice paths, and red points represent ending points of them.}
\end{figure}

Note that the starting points for our paths correspond to the unit holes indexed by \(P\), as well as to the unit holes indexed by elements of \(Q\) that do not have a corresponding unit hole indexed by an element of \(R\) next to them. Similarly, the ending points correspond to the unit holes indexed by \(S\), and to the unit holes indexed by elements of \(R\) that do not have a corresponding unit hole indexed by an element of \(Q\) next to them.

Choose a system of coordinates as indicated in Figure 3.2 (we place \(L_{n,x}(P,Q,R,S)\) so that the inner unit holes are along the first bisector, and its base is half a lattice spacing above the \(x\)-axis; \(\overline{L}_{n,x}(P,Q,R,S)\) is placed so that the inner holes are along the first bisector and its left side is half a lattice spacing to the right of the \(y\)-axis).

Then by the above mentioned interpretation, we obtain that the number of lozenge tilings of \(L_{n,x}(P,Q,R,S)\) is equal to the number of families of non-intersecting lattice paths on \(\mathbb{Z}^2\) with starting points \(\{(x+1,x+p)|p \in P\}\cup \{(x+q,x+q)|q \in Q\cap R'\}\) and ending points \(\{(x+r,x+r)|r \in Q'\cap R\}\cup \{(x+s,1)|s \in S\})\), allowed to take steps south or east.

Among the \(|P|\) lattice paths starting at \(\{(x+1,x+p)|p \in P\}\), \(|Q'\cap R|\) of them end at points in the set \(\{(x+r,x+r)|r \in Q'\cap R\}\), and the remaining \(|P|-|Q'\cap R|\) pass through the first bisector\footnote{ It readily follows from the interpretation of tilings as paths of lozenges that if a tiling of \(L_{n,x}(P,Q,R,S)\) exists (as we are assuming), then \(|P|-|Q'\cap R|\geq0\).}.

We now partition the collection of families of non-intersecting lattice paths corresponding to the tilings of our regions \(L_{n,x}(P,Q,R,S)\) (we will soon do the same for the regions \(\overline{L}_{n,x}(P,Q,R,S)\)) according to the set of points on the first bisector that they pass through. Note that  
this set of points has coordinates of the form \(\{(x+u,x+u)|u\in U\}\), where \(U\subset [n]\setminus(Q\cup R)\) and \(|U|=|P|-|Q'\cap R|\).

Each family of non-intersecting lattice paths that pass through \(\{(x+u,x+u)|u\in U\}\) can be thought of as a pair of families non-intersecting lattice paths: One starting from \(\{(x+1,x+p)|p \in P\}\) and ending at \(\{(x+r,x+r)|r \in (Q'\cap R)\cup U\}\), and the other starting from \(\{(x+q,x+q)|q \in (Q\cap R')\cup U\}\) and ending at \(\{(x+s,1)|s \in S\})\) (see the left picture in Figure 3.2). 
Note that the starting and ending points of both these families of lattice paths satisfy the compatibility condition that Lindstr\"{o}m-Gessel-Viennot theorem requires.

Let \(p_1, p_2,\dotsc, p_k\) be the elements of \(P\), \(r^U_1, r^U_2,\dotsc, r^U_k\) the elements of \((Q'\cap R)\cup U\), \(q^U_1, q^U_2,\dotsc, q^U_l\) the elements of \((Q\cap R')\cup U\) and \(s_1, s_2,\dotsc, s_l\) the elements of \(S\), where elements are written in increasing order.

Then, by Lindstr\"{o}m-Gessel-Viennot theorem, we can express the number of lozenge tilings of \(L_{n,x}(P,Q,R,S)\) as follows\footnote{We are also using the fact that for \(a,b,c,d\in \mathbb{Z}\), the number of lattice paths on \(\mathbb{Z}^2\) from \((a,b)\) to \((c,d)\) taking steps south or east is 
\begin{equation*}
\begin{cases}
   \binom{(c-a)+(b-d)}{c-a} & \text{if \((c-a) + (b-d) \geq 0\) \text{and} \((c-a) \geq 0\)}, \\
   \binom{(c-a)+(b-d)}{b-d} & \text{if \((c-a) + (b-d) \geq 0\) \text{and} \((b-d) \geq 0\)}.
\end{cases}    
\end{equation*}
}.

\begin{equation}
\begin{aligned}
    &M(L_{n,x}(P,Q,R,S))\\
    &\begin{aligned}
                    =\sum_{U} &\text{det}\Bigg[\binom{(x+r^U_{j})-(x+1)+(x+p_{i})-(x+r^U_{j})}{(x+r^U_{j})-(x+1)}\Bigg]\\
                    &\ \text{det}\Bigg[\binom{(x+s_{j})-(x+q^U_{i})+(x+q^U_{i})-1}{(x+q^U_{i})-1}\Bigg]\\
                \end{aligned}\\
    &=\sum_{U} \text{det}\Bigg[\binom{p_{i}-1}{r^U_{j}-1}\Bigg] \ \text{det}\Bigg[\binom{x+s_{j}-1}{x+q^U_{i}-1}\Bigg]
\end{aligned}
\end{equation}
where the summation is over all sets \(U\subset [n]\setminus(Q\cup R)\) with \(|U|=|P|-|Q'\cap R|\).

Similarly, \(M(\overline{L}_{n,x}(P,Q,R,S))\) is the same as the number of non-intersecting lattice paths on \(\mathbb{Z}^2\) from starting points \(\{(1,x+p)|p \in P\}\cup \{(x+q,x+q)|q \in Q\cap R'\}\) to ending points \(\{(x+r,x+r)|r \in Q'\cap R\}\cup \{(x+s,x+1)|s \in S\})\), where paths can only move  south or east. These paths are also in bijection with pairs of families non-intersecting lattice paths (see the right picture in Figure 3.2).

By the same partitioning and the Lindstr\"{o}m-Gessel-Viennot theorem, we have
\begin{equation}
\begin{aligned}
    &M(\overline{L}_{n,x}(P,Q,R,S))\\
    &=\sum_{U} \text{det}\Bigg[\binom{x+p_{i}-1}{x+r^U_{j}-1}\Bigg] \ \text{det}\Bigg[\binom{s_{j}-1}{q^U_{i}-1}\Bigg]
\end{aligned}
\end{equation}
where the summation is over the same sets \(U\) as in (3.2).

For such \(U\), by the definition of binomial coefficients and by the linearity of the determinant in rows and columns, we have 
\begin{equation}
\begin{aligned}
    \text{det}\Bigg[\binom{x+p_{i}-1}{x+r^U_{j}-1}\Bigg]=\text{det}\Bigg[\frac{(p_{i})_{x}}{(r^U_{j})_{x}} \ &\binom{p_{i}-1}{r^U_{j}-1}\Bigg]=\frac{\displaystyle \prod_{i=1}^{k} (p_{i})_{x}}{\displaystyle \prod_{j=1}^{k} (r^U_{j})_{x}} \ \text{det}\Bigg[\binom{p_{i}-1}{r^U_{j}-1}\Bigg]\\
    &=\frac{\displaystyle \prod_{p\in P} (p)_{x}}{\displaystyle \prod_{r\in (Q'\cap R)\cup U} (r)_{x}} \ \text{det}\Bigg[\binom{p_{i}-1}{r^U_{j}-1}\Bigg]
\end{aligned}
\end{equation}
and
\begin{equation}
\begin{aligned}
    \text{det}\Bigg[\binom{x+s_{j}-1}{x+q^U_{i}-1}\Bigg]=\text{det}\Bigg[\frac{(s_{j})_{x}}{(q^U_{i})_{x}}\ &\binom{s_{j}-1}{q^U_{i}-1}\Bigg]=\frac{\displaystyle \prod_{j=1}^{l} (s_{j})_{x}}{\displaystyle \prod_{i=1}^{l} (q^U_{i})_{x}}\ \text{det}\Bigg[\binom{s_{j}-1}{q^U_{i}-1}\Bigg]\\
    &=\frac{\displaystyle \prod_{s\in S} (s)_{x}}{\displaystyle \prod_{q\in (Q\cap R')\cup U} (q)_{x}}\ \text{det}\Bigg[\binom{s_{j}-1}{q^U_{i}-1}\Bigg].
\end{aligned}
\end{equation}

Since we factored out non-zero factors in (3.4) and (3.5), we have

\begin{equation}
\begin{aligned}
    &\text{det}\Bigg[\binom{p_{i}-1}{r^U_{j}-1}\Bigg] \ \text{det}\Bigg[\binom{x+s_{j}-1}{x+q^U_{i}-1}\Bigg]\neq 0\\
    &\iff \text{det}\Bigg[\binom{x+p_{i}-1}{x+r^U_{j}-1}\Bigg] \ \text{det}\Bigg[\binom{s_{j}-1}{q^U_{i}-1}\Bigg]\neq 0.
\end{aligned}    
\end{equation}

By (3.2), existence of a lozenge tiling of the region \(L_{n,x}(P,Q,R,S)\) guarantees existence of a set \(U\) that satisfies the condition in (3.6). For any such \(U\), by (3.4) and (3.5), we have
\begin{equation}
\begin{aligned}
    \frac{\text{det}\Big[\binom{x+p_{i}-1}{x+r^U_{j}-1}\Big] \ \text{det}\Big[\binom{s_{j}-1}{q^U_{i}-1}\Big]}{\text{det}\Big[\binom{p_{i}-1}{r^U_{j}-1}\Big] \ \text{det}\Big[\binom{x+s_{j}-1}{x+q^U_{i}-1}\Big]}
    &=\frac{\displaystyle \prod_{p\in P} (p)_{x}}{\displaystyle \prod_{r\in (Q'\cap R)\cup U} (r)_{x}} \ \frac{\displaystyle \prod_{q\in (Q\cap R')\cup U} (q)_{x}}{\displaystyle \prod_{s\in S} (s)_{x}}\\
    &=\frac{\displaystyle \prod_{p\in P} (p)_{x} \ \prod_{q\in Q} (q)_{x}}{\displaystyle \prod_{r\in R} (r)_{x}\ \prod_{s\in S} (s)_{x}}.
\end{aligned}    
\end{equation}

Since the ratio does not depend on \(U\), by (3.2), (3.3) and (3.7), we have
\begin{equation}
    \frac{M(\overline{L}_{n,x}(P,Q,R,S))}{M(L_{n,x}(P,Q,R,S))}=\frac{\displaystyle \prod_{p\in P} (p)_{x} \ \prod_{q\in Q} (q)_{x}}{\displaystyle \prod_{r\in R} (r)_{x} \ \prod_{s\in S} (s)_{x}}.
\end{equation}

This completes the proof. \(\square\)

\medskip
From (3.2), (3.3) and (3.6), and using also the fact that the factors in the summands of (3.2) and (3.3) are non-negative (as, by the Lindstr\"{o}m-Gessel-Viennot theorem, each represents the cardinality of a certain family of non-intersecting lattice paths), it follows that
\begin{equation}
    M(L_{n,x}(P,Q,R,S))\neq 0 \iff M(\overline{L}_{n,x}(P,Q,R,S))\neq 0.
\end{equation}

The proof of the main theorem is almost the same as that of Lemma 2.1. We partition the set of lozenge tilings of both the \(H\)- and \(\overline{H}\)-regions according to the positions of the lozenges that cross one of three line segments (we specify in the next paragraph what these are). Again, we group partition classes from two families in pairs. By using Lemma 3.1, we then deduce that the ratio of cardinalities of corresponding partition classes is equal to the right-hand side of (2.1).
\medskip

\textit{Proof of Theorem 2.1.} We partition the set of lozenge tilings of \(H_{n,x}(\textbf{A},\textbf{B})\) according to positions of lozenges crossing the following three line segments: The one connecting \(R\) and the right vertex of the hexagon\footnote{ Recall that $U$, $L$ and $R$ are the vertices of the central triangle (see Figure 3.3).}, the one connecting \(L\) and the bottom left vertex of the hexagon, and the one connecting \(U\) and top left vertex of the hexagon. Let \(W_2\), \(W_4\) and \(W_6\) be the sets of indices of unit segments crossed by lozenges on each of these segments, respectively.
For \(\textbf{W}\coloneqq(W_2, W_4, W_6)\), let \(\textbf{A}_\textbf{W} \coloneqq (A_1\cup W_6, A_2, A_3\cup W_2, A_4, A_5\cup W_4, A_6)\) and \(\textbf{B}^\textbf{W} \coloneqq (B_1, B_2\cup W_2, B_3, B_4\cup W_4, B_5, B_6\cup W_6)\) (see Figure 2.3 to recall what the $A_i$’s and $B_i$’s are). Then\footnote{Here, we are also using that the set of lozenge tilings of a given region with some fixed lozenges at specified positions is clearly in bijection with the set of lozenge tilings of the same region with additional lozenge-shaped holes at those positions.}

\begin{figure}
    \centering
    \includegraphics[width=10.5cm]{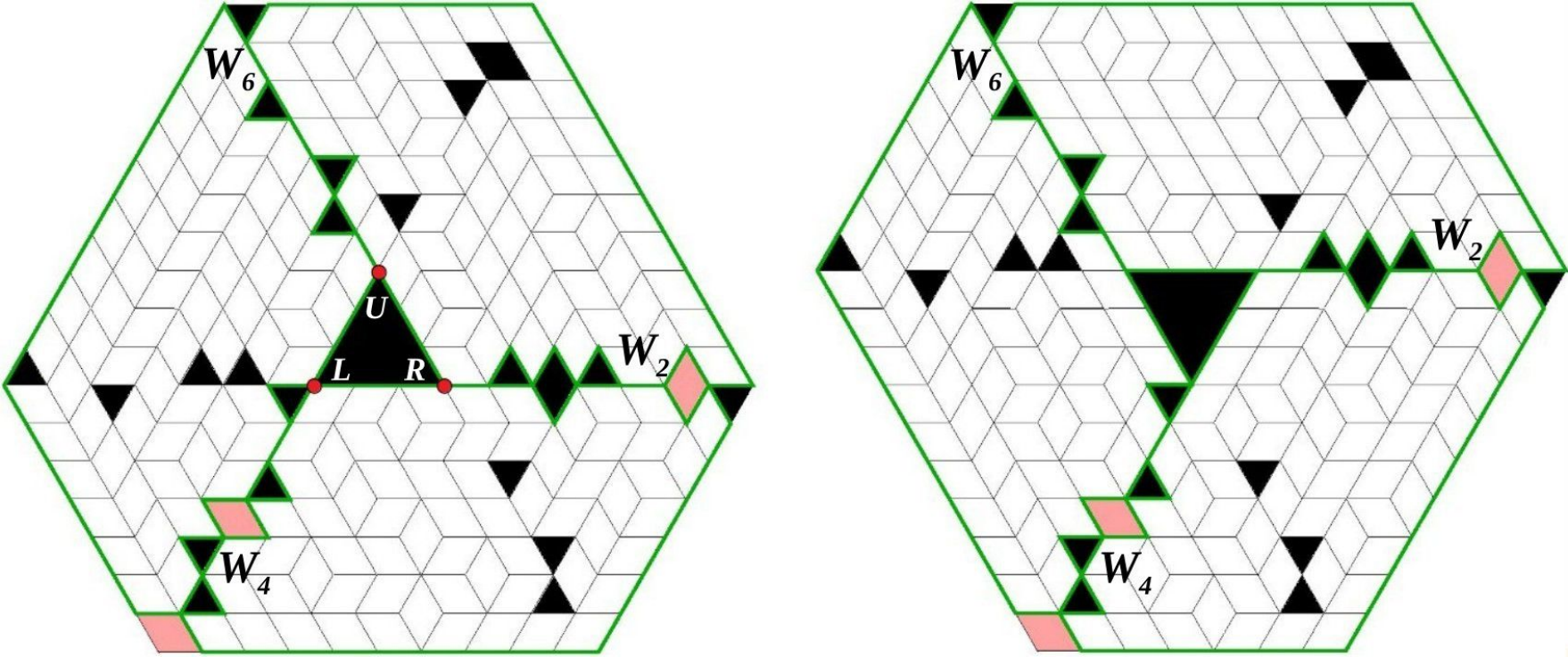}
    \caption{The regions \(H_{7,3}(\textbf{A}_\textbf{W},\textbf{B}^\textbf{W})\) (left) and \(\overline{H}_{7,3}(\textbf{A}_\textbf{W},\textbf{B}^\textbf{W})\) (right), with \(\textbf{W}=(W_2, W_4, W_6)=(\{6\},\{4,7\}, \emptyset)\), and sample lozenge tilings. \(H_{7,3}(\textbf{A}_\textbf{W},\textbf{B}^\textbf{W})\) is split into three \(L\)-regions, and \(\overline{H}_{7,3}(\textbf{A}_\textbf{W},\textbf{B}^\textbf{W})\) is split into three \(\overline{L}\)-regions. The shaded lozenges are the ones specified by \(\textbf{W}\).}
\end{figure}

\begin{equation}
    M(H_{n,x}(\textbf{A},\textbf{B}))=\sum_{\textbf{W}} M(H_{n,x}(\textbf{A}_\textbf{W},\textbf{B}^\textbf{W}))
\end{equation}
where sum runs over all triples \(\textbf{W}\) for which \(W_{2i}\cap (A_{2i+1}\cup B_{2i})=\emptyset\) for \(i=1,2,3\) (where \(A_7 := A_1\)) and \(M(H_{n,x}(\textbf{A}_\textbf{W},\textbf{B}^\textbf{W}))>0\).

Observe that for any such \(\textbf{W}\), the region \(H_{n,x}(\textbf{A}_\textbf{W},\textbf{B}^\textbf{W})\) can be split into three subregions — separated by the three line segments  mentioned above (see the left picture in Figure 3.3) — which are \(L\)-regions of the kind we dealt with in the previous lemma. Furthermore, because of our partitioning, for any tiling of \(H_{n,x}(\textbf{A}_\textbf{W},\textbf{B}^\textbf{W})\), each of these regions is tiled internally.
Hence, if we set \(W_0 := W_6\), we have
\begin{equation}
    M(H_{n,x}(\textbf{A},\textbf{B}))=\displaystyle \sum_{\textbf{W}} \Bigg[\prod_{i=1}^{3} M(L_{n,x}(A_{2i-1}\cup W_{2i-2}, B_{2i-1}, A_{2i}, B_{2i}\cup W_{2i}))\Bigg].
\end{equation}

By the same argument (see right picture in Figure 3.3),
\begin{equation}
    M(\overline{H}_{n,x}(\textbf{A},\textbf{B}))=\displaystyle \sum_{\textbf{W}} \Bigg[\prod_{i=1}^{3} M(\overline{L}_{n,x}(A_{2i-1}\cup W_{2i-2}, B_{2i-1}, A_{2i}, B_{2i}\cup W_{2i}))\Bigg]
\end{equation}
where the sum is over all triples \(\textbf{W}=(W_2, W_4, W_6)\subset [n]^{3}\) such that \(W_{2i}\cap (A_{2i+1}\cup B_{2i})=\emptyset\) for \(i=1,2,3\) and \(M(\overline{H}_{n,x}(\textbf{A}_\textbf{W},\textbf{B}^\textbf{W}))>0\).

\medskip
In fact, the summation ranges for \(\textbf{W}\) in (3.11) and (3.12) turn out to be the same. Indeed,
\begin{equation}
\begin{aligned}
    &M(H_{n,x}(\textbf{A}_\textbf{W},\textbf{B}^\textbf{W}))\neq 0\\
    &\iff \displaystyle M(L_{n,x}(A_{2i-1}\cup W_{2i-2}, B_{2i-1}, A_{2i}, B_{2i}\cup W_{2i}))\neq 0, \forall i=1,2,3\\
    &\iff \displaystyle M(\overline{L}_{n,x}(A_{2i-1}\cup W_{2i-2}, B_{2i-1}, A_{2i}, B_{2i}\cup W_{2i}))\neq 0, \forall i=1,2,3\\
    &\iff M(\overline{H}_{n,x}(\textbf{A}_\textbf{W},\textbf{B}^\textbf{W}))\neq 0,
\end{aligned}    
\end{equation}
where at the second step we used (3.9).

\medskip
Also, by (3.10), the existence of a lozenge tiling of the region \(H_{n,x}(\textbf{A},\textbf{B})\) guarantees the existence of a \(\textbf{W}\) that satisfies the condition in (3.13). For any such \(\textbf{W}\), we analyze the ratio of corresponding summands in (3.11) and (3.12). By Lemma 3.1, we have 

\begin{equation}
\begin{aligned}
    &\frac{\prod_{i=1}^{3} M(\overline{L}_{n,x}(A_{2i-1}\cup W_{2i-2}, B_{2i-1}, A_{2i}, B_{2i}\cup W_{2i}))}{\prod_{i=1}^{3} M(L_{n,x}(A_{2i-1}\cup W_{2i-2}, B_{2i-1}, A_{2i}, B_{2i}\cup W_{2i}))}\\
    &=\displaystyle \prod_{i=1}^{3} \Bigg[\frac{\displaystyle \prod_{a\in A_{2i-1}\cup W_{2i-2}} (a)_{x} \ \prod_{b\in B_{2i-1}} (b)_{x}}{\displaystyle \prod_{a\in A_{2i}} (a)_{x} \ \prod_{b\in B_{2i}\cup W_{2i}} (b)_{x}}\Bigg]\\
    &=\displaystyle \prod_{i=1}^{3} \Bigg[\displaystyle \frac{\prod_{a\in A_{2i-1}} (a)_{x} \ \prod_{b\in B_{2i-1}} (b)_{x}}{\prod_{a\in A_{2i}} (a)_{x} \ \prod_{b\in B_{2i}} (b)_{x}}\Bigg]\\    
    &=\frac{\displaystyle \prod_{a\in A_{o}}(a)_x \ \displaystyle \prod_{b\in B_{o}}(b)_x}{\displaystyle \prod_{a\in A_{e}}(a)_x \ \displaystyle \prod_{b\in B_{e}}(b)_x},
\end{aligned}
\end{equation}
where at the second to last equality we used \(W_{2i}\cap(A_{2i+1}\cup B_{2i})=\emptyset\), so that factors corresponding to elements of \(W_{2}\cup W_{4} \cup W_{6}\) cancel out completely.

Note that the right hand side of (3.14) does not depend on \(\textbf{W}\). Hence, from (3.11), (3.12) and (3.14), we have
\begin{equation}
    \frac{M(\overline{H}_{n,x}(\textbf{A},\textbf{B}))}{M(H_{n,x}(\textbf{A},\textbf{B}))}=\frac{\displaystyle \prod_{a\in A_{o}}(a)_x \ \displaystyle \prod_{b\in B_{o}}(b)_x}{\displaystyle \prod_{a\in A_{e}}(a)_x \ \displaystyle \prod_{b\in B_{e}}(b)_x}.
\end{equation}

This completes the proof of the first part.

\smallskip
For the second part, recall that a cyclically symmetric lozenge tiling is a tiling that is invariant under rotation by 120\(^{\circ}\).
Hence, the position of lozenges that cross the three lines are also invariant under this rotation (i.e. \(W_2=W_4=W_6\)).

\begin{figure}
    \centering
    \includegraphics[width=10.5cm]{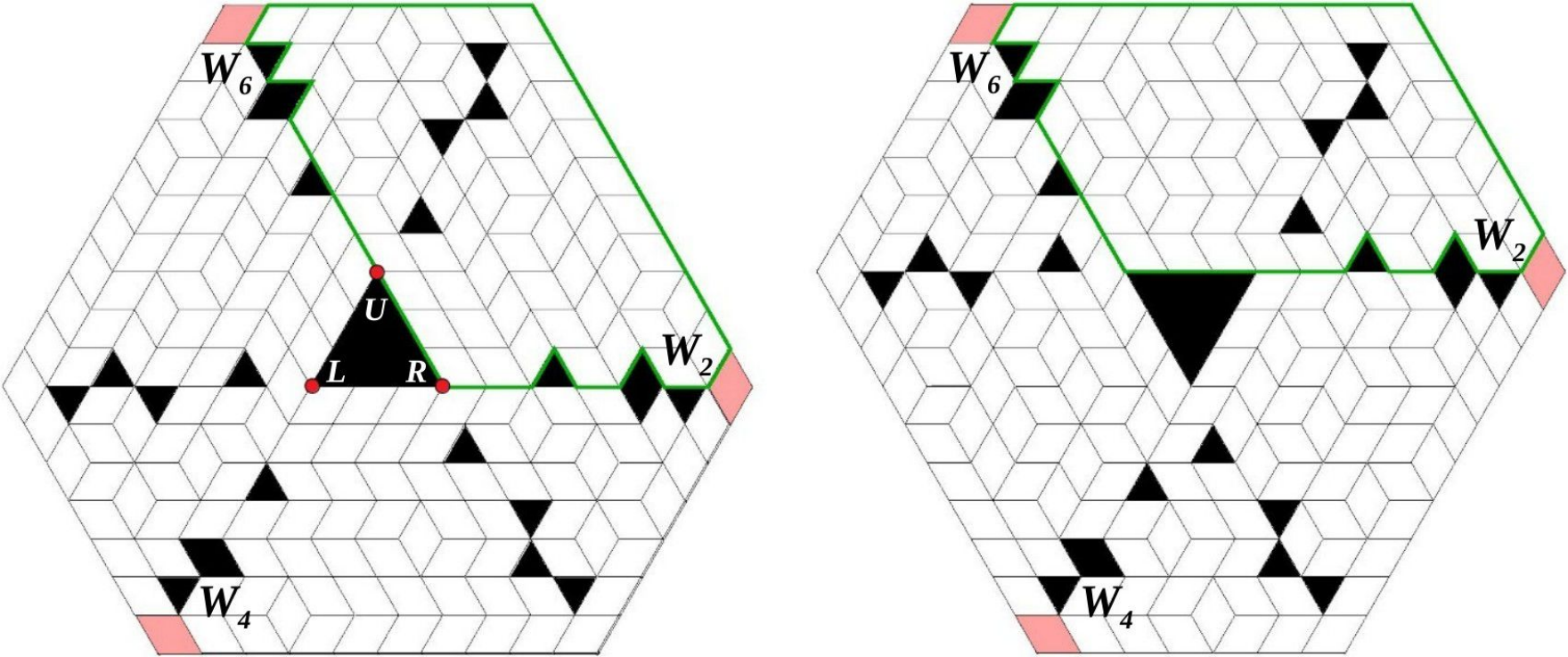}
    \caption{The cyclically symmetric regions \(H_{7,3}(\textbf{A}_\textbf{W},\textbf{B}^\textbf{W})\) (left) and \(\overline{H}_{7,3}(\textbf{A}_\textbf{W},\textbf{B}^\textbf{W})\) (right), with \(\textbf{W}=(W_2, W_4, W_6)=(\{7\},\{7\},\{7\})\) and examples of their cyclically symmetric lozenge tilings. Cyclically symmetric tilings are determined by tilings of the indicated subregions.}
\end{figure}

One can readily see that cyclically symmetric tilings of \(H_{n,x}(\textbf{A}_\textbf{W},\textbf{B}^\textbf{W})\) are in bijection with tilings of \(L_{n,x}(A_{1}\cup W_{2}, B_{1}, A_{2}, B_{2}\cup W_{2})\) (see Figure 3.4). Thus, by the same partitioning that we used in the first part, we have
\begin{equation}
    M_{r}(H_{n,x}(\textbf{A},\textbf{B}))=\displaystyle \sum_{W_2} M(L_{n,x}(A_{1}\cup W_{2}, B_{1}, A_{2}, B_{2}\cup W_{2})).
\end{equation}
where the summation is over the all \(W_2\) such that \(W_{2}\cap (A_{1}\cup B_{2})=\emptyset\).

By the same argument,
\begin{equation}
    M_{r}(\overline{H}_{n,x}(\textbf{A},\textbf{B}))=\displaystyle \sum_{W_2} M(\overline{L}_{n,x}(A_{1}\cup W_{2}, B_{1}, A_{2}, B_{2}\cup W_{2})).
\end{equation}

As in the proof of the first part, the summation in (3.17) is over the same \(W_2\) as in (3.16). Furthermore, the existence of a cyclically symmetric lozenge tiling of the region \(H_{n,x}(\textbf{A},\textbf{B})\) guarantees that the right-hand side of (3.16) contains a non-zero summand.
Then, by the same manipulation we did in the proof of the first part, we have
\begin{equation}
    \frac{M_{r}(\overline{H}_{n,x}(\textbf{A},\textbf{B}))}{M_{r}(H_{n,x}(\textbf{A},\textbf{B}))}=\displaystyle \frac{\prod_{a\in A_{1}} (a)_{x} \ \prod_{b\in B_{1}} (b)_{x}}{\prod_{a\in A_{2}} (a)_{x} \ \prod_{b\in B_{2}} (b)_{x}}.
\end{equation}

The remaining equality in (2.2) follows directly from the first part. \qedsymbol\\

\textbf{Remark.} In this paper, we showed that the ratio between the numbers of lozenge tilings of the two regions \(H_{n,x}(\textbf{A},\textbf{B})\) and \(\overline{H}_{n,x}(\textbf{A},\textbf{B})\) is expressed as a simple product formula. We also provided a corresponding formula for the cyclically symmetric tilings, which turns out to be the cube root of the former formula. It is then natural to ask what happens for other symmetry classes. Two more symmetry classes make sense in our context: vertically symmetric tilings, and cyclically symmetric and vertically symmetric tilings.

It turns out that the same arguments we used in the proof of Theorem 2.1, together with a result due independently to Condon [14] and Lai [26, Theorem 1.3] allow us to prove the following formulas.

\smallskip
If \(H_{n,x}(\textbf{A},\textbf{B})\) is vertically symmetric\footnote{ I.e. it is invariant under reflection across a vertical line.} and has a vertically symmetric lozenge tiling (i.e. a tiling invariant under reflection across a vertical line), then 
\begin{equation}
    \frac{M_{|}(\overline{H}_{n,x}(\textbf{A},\textbf{B}))}{M_{|}(H_{n,x}(\textbf{A},\textbf{B}))}=\frac{\displaystyle \prod_{a\in A_{o}}(a)_x}{\displaystyle \prod_{a\in A_{e}}(a)_x}=\sqrt{\frac{M(\overline{H}_{n,x}(\textbf{A},\textbf{B}))}{M(H_{n,x}(\textbf{A},\textbf{B}))}},
\end{equation}
where $M_|(R)$ denotes the number of vertically symmetric tilings of the region $R$.

Similarly, if \(H_{n,x}(\textbf{A},\textbf{B})\) is cyclically symmetric and vertically symmetric, and has a cyclically symmetric and vertically symmetric lozenge tiling, then, denoting by $M_{r,|}(R)$ the number of tilings of the region $R$ that are invariant under both reflection across a vertical line and rotation by 120\(^{\circ}\), we have
\begin{equation}
\begin{aligned}
    \frac{M_{r,|}(\overline{H}_{n,x}(\textbf{A},\textbf{B}))}{M_{r,|}(H_{n,x}(\textbf{A},\textbf{B}))}=\frac{\displaystyle \prod_{a\in A_{1}}(a)_x}{\displaystyle \prod_{a\in A_{2}}(a)_x}&=\sqrt{\frac{M_{r}(\overline{H}_{n,x}(\textbf{A},\textbf{B}))}{M_{r}(H_{n,x}(\textbf{A},\textbf{B}))}}\\
    &=\sqrt[3]{\frac{M_{|}(\overline{H}_{n,x}(\textbf{A},\textbf{B}))}{M_{|}(H_{n,x}(\textbf{A},\textbf{B}))}}\\
    &=\sqrt[6]{\frac{M(\overline{H}_{n,x}(\textbf{A},\textbf{B}))}{M(H_{n,x}(\textbf{A},\textbf{B}))}}.
\end{aligned}
\end{equation}

Just as from the second part of Theorem 2.1 one could get a simple proof for the enumeration of cyclically symmetric lozenge tilings of hexagons with a shamrock removed from the center (see the end of Section 2), the above two identities can be used to deduce simple proofs for the formulas enumerating vertically symmetric lozenge tilings and cyclically symmetric and vertically symmetric lozenge tilings of those regions, which were first proved by Ciucu in [6] and [5], respectively.

As one of the referees pointed out, the simplicity of the formulas in our results makes one wonder if there are even simpler proofs for them, for instance in the style of a bijective proof. It would be interesting if such proofs could be found, but this seems to be a hard goal to achieve. 

\medskip

\textbf{Acknowledgments.}
The author thanks his advisor, Professor Mihai Ciucu, for useful discussions and his encouragement. This paper could not have been written without his guidance and discussions with him about the bowtie squeezing theorem. The author also thanks anonymous reviewers for carefully reading the original version of the paper and giving helpful comments. Also, the author expresses gratitude to Jeff Taylor for frequent and helpful assistance.
Lastly, the author thanks David Wilson because his program, \textit{Vaxmacs}, helped the author a lot when he made the observation that led to Theorem 2.1.

\end{document}